\documentclass{article}
\usepackage{amsmath}
\usepackage{amssymb}
\usepackage{amsthm}
\usepackage{graphicx}
\begin{document}
\newtheorem{theorem}{Theorem}
\newtheorem{prop}{Proposition}
\newtheorem{proposition}{Proposition}
\newtheorem{lemma}{Lemma}
\newtheorem*{fact}{Fact}
\theoremstyle{definition}
\newtheorem{definition}{Definition}
\newtheorem{axiom}{Axiom}
\newtheorem{exercise}{Exercise}
\bibliographystyle{amsplain}
\title{An Introduction to Smooth Infinitesimal Analysis}
\author{Michael O'Connor}
\maketitle
\newcommand{\N}{\mathbb{N}}
\tableofcontents
\section{Why I Wrote This}
The primary reason that I wrote this was to have a freely available version of this material on the web.  Nothing
in this article is due to me (except for any mistakes).  My primary sources were~\cite{primer}~and~\cite{moerdijk};
most everything up to Stokes's theorem is from~\cite{primer} and most everything after and including Stokes's theorem is from~\cite{moerdijk}. I would also
direct the reader to~\cite{kock}~and~\cite{shulman}.

\section{Introduction and Motivation}
Many mathematicians, from Archimedes to Leibniz to Euler and beyond, made use of infinitesimals in their arguments.  These were later replaced rigorously with limits, but many people still
find it useful to think and derive with infinitesimals.

Unfortunately, in most informal setups the existence of infinitesimals is technically contradictory, so it can
be difficult to grasp the means by which one fruitfully manipulates them.  It would be useful to have an axiomatic framework with the following properties:  

1. It is consistent.  

2. The system acts as a good ``intuition pump''
for the real world.  In particular, this entails that if you prove something in the system, then while it won't necessarily be \emph{true} in the real world, there should be a high probability that it's \emph{morally true} in the real world, i.e., with some extra assumptions it becomes true.  It should also ideally entail that
many of the proofs of Archimedes, et al., involving infinitesimals can be formulated as is (or close to ``as is'').

``Smooth infinitesimal analysis'' is one attempt to satisfy these conditions.

\section{Axioms and Logic}

Consider the following axioms:

\begin{axiom}$R$ is a set, 0 and 1 are elements of $R$ and $+$ and $\cdot$ are binary operations on $R$.  The structure $\langle R, +, \cdot, 0, 1\rangle$ is a commutative ring with unit.  

Furthermore, we have that $\forall x\, ((x \ne 0) \implies (\exists y\, xy = 1))$, but I don't want to call
$R$ a field for a reason I'll discuss in a moment.
\end{axiom}
\begin{axiom} There is a transitive irreflexive relation $<$ on $R$.  It satisfies $0 < 1$, and for all $x$, $y$, and
$z$, we have $x < y \implies x + z < y + z$ and ($x < y$ and $z > 0)\implies xz < yz$.

It also satisfies $\forall x, y\, (x\ne y) \implies (x > y \vee x < y)$,
but I don't want to call $<$ total, for a reason I'll discuss in a moment.
\end{axiom}
\begin{axiom} For all $x > 0$ there is a unique $y > 0$ such that $y^2 = x$.
\end{axiom}
\begin{axiom}[Kock-Lawvere Axiom] Let $D = \{d \in R\mid d^2 = 0\}$.  Then for all functions $f$ from $D$ to $R$,
and all $d\in D$, there is a unique $a\in R$ such that $f(d) = f(0) + d\cdot a$.
\end{axiom}

After reading the Kock-Lawvere Axiom you are probably quite puzzled.  In the first place, we
can easily prove that $D = \{0\}$: Let $d\in D$. For a proof by contradiction, assume that $d\ne 0$, then there is a $d^{-1}$ and if $d^2$ equalled 0, we would have $d = d^2 d^{-1} = 0$.

For an alternate proof that $D = \{0\}$: Again assume that $d\ne 0$ for a contradiction. Then $d > 0$
or $d < 0$.  In the first case, $d^2 > 0$, so $d\ne 0$ (since $<$ is irreflexive). In the second case,
we have $0 < -d$ by adding $-d$ to both sides, and again $d^2 > 0$.

Now, if $D = \{0\}$, then for \emph{any} $a\in R$, and any function $f$ from $D$ to $R$, we have
$f(d) = f(0) + d\cdot a$ for all $d\in D$.  This contradicts the uniqueness of $a$.  Therefore, the axioms
presented so far are contradictory.

However, we have the following surprising fact.

\begin{fact} There is a form of set theory (called a \emph{local set theory}, or \emph{topos logic}) which has
its underlying logic restricted (to a logic called \emph{intuitionistic logic}) under which Axioms 1 through 4
(and also the axioms to be presented later in this paper)
taken together are
consistent.
\end{fact}

\begin{definition}[Smooth Infinitesimal Analysis] Smooth Infinitesimal Analysis (SIA) is the system whose axioms are
those sentences marked as Axioms in this paper and whose logic is that alluded to in the above theorem.
\end{definition}

References for this theorem are~\cite{moerdijk}~and~\cite{kock}.  References for topos logic specifically are~\cite{localset} and~\cite{moerdijkmaclane}.

Essentially, intuitionistic logic disallows proof by contradiction (which was used in both proofs that $D = \{0\}$ above) and its equivalent brother, the law of the excluded middle, which says that for any proposition $P$, $P\vee \neg P$ holds.

I won't formally define intuitionistic logic or topos logic here as it would take too much space
and there's no real way to understand it except by seeing examples anyway.
If you avoid proofs by contradiction and 
proofs using the law of the excluded middle (which usually come up in ways like: ``Let $x\in R$.  Then either
$x = 0$ or $x\ne 0$.\ldots''), you will be okay.

But before we go further we might ask, ``what does this logic have to do with the real world anyway?''  Possibly nothing, but recall that our goals above do not require that we work with ``real'' objects; just that we have a consistent system which will act as a good ``intuition pump'' about the real world.  We are guaranteed that the system is consistent by a theorem; for the second condition each person will have to judge for themselves.

To conclude this section, it should now be clear why I didn't want to call $R$ a field and $<$ a total order:
Even though we have $\forall x\, ((x\ne 0)\implies \text{$x$ invertible})$, we can't conclude from that that
$\forall x\, ((x = 0) \vee (\text{$x$ invertible}))$, because the proof of the latter from the former uses
the law of the excluded middle.  Calling $R$ a field would unduly give the impression that the latter is true.

For the rest of this paper I will generally work within SIA (except, obiviously, when I announce new
axioms or make remarks about SIA).

\section{Single-Variable Calculus}

\subsection{An Important Lemma}
This lemma is easy to prove, but because it is used over and over again, I'll isolate it here:
\begin{lemma}[Microcancellation] Let $a, b\in R$. If for all $d\in D$ we have $ad = bd$, then $a = b$.
\end{lemma}
\begin{proof}
Let $f\in R^D$ be given by $f(d) = ad = bd$.  Then by the uniqueness condition of the Kock-Lawvere axiom,
we have that $a = b$.
\end{proof}

\subsection{Basic Rules}
Let $f$ be a function from $R$ to $R$, and let $x\in R$.  We may defined a function $g$ from $D$ to $R$ as
follows: for all $d\in D$, let $g(d) = f(x + d)$.  Then the Kock-Lawvere axiom tells us that there is a unique
$a$ so that $g(d) = g(0) + ad$ for all $d\in D$.  Thus, we have that for all functions $f$ from $R$ to $R$
and all $x\in R$, there is a unique $a$ so that $f(x + d) = f(x) + ad$ for all $d$.  We define $f'(x)$ to be
this $a$.

We thus have the following fundamental fact: 

\begin{proposition}[Fundamental Fact about Derivatives]
For all $f\in R^R$, all $x\in R$, and all $d\in D$,
\[f(x + d) = f(x) + f'(x)d\] and furthermore, $f'(x)$ is the unique real number with that property.
\end{proposition}

\begin{proposition} Let $f$, $g\in R^R$, $c\in R$. Then:

1. $(f + g)' = f' + g'$

2. $(cf)' = cf'$

3. $(fg)' = f'g + fg'$.

4. If for all $x$, $g(x) \ne 0$, then $(f/g)' = (gf' - fg')/g^2$.

5. $(f\circ g)' = (f'\circ g)\cdot g'$.
\end{proposition}
\begin{proof}
I'll prove 3 and 5 and leave the rest as exercises. 

To prove 3: Let $x\in R$ and $d\in D$.  Let $h(x) = f(x)g(x)$. Then 
\begin{displaymath}
h(x + d) = f(x + d)g(x + d) = (f(x) + f'(x)d)(g(x) + g'(x)d)
\end{displaymath} 
which, multiplying out and using $d^2 = 0$, is equal to 
\[f(x)g(x) + d(f'(x)g(x) + f(x)g'(x)) = h(x) +d(f'(x)g(x) + f(x)g'(x)).\]  On the other hand, we know that
$h(x + d) = h(x) + h'(x)d$, so \[h'(x)d = d(f'(x)g(x) + f(x)g'(x)).\]  Since $d$ was an arbitrary element of $D$, we may use
microcancellation, and we obtain $h'(x) = f'(x)g(x) + f(x)g'(x)$.

To prove 5: Let $x\in R$ and $d\in D$. Then \[f(g(x + d)) = f(g(x) + g'(x)d).\]  Now, $g'(x)d$ is in
$D$ (since $(g'(x)d)^2 = d^2(g'(x))^2 = 0$), so \[f(g(x) + g'(x)d) = f(g(x)) + g'(x)f'(g(x))d.\]
As before, this gives us that 
 $g'(x)f'(g(x))$ is the derivative of $f(g(x))$.
\end{proof}

In order to do integration, let's add the following axiom:

\begin{axiom} For all $f\in R^R$ there is a unique $g\in R^R$ such that $g' = f$ and $g(0) = 0$.  We write
$g(x)$ as $\int_0^x f(t)\,dt$.
\end{axiom}

We can now derive the rules of integration in the usual way by inverting the rules of differentiation.

\subsection{Deriving formulas for Arclength, etc.}

I'd now like to derive the formula for the arclength of the graph of a function $y = f(x)$ (say,
from $x = 0$ to $x = 1$).  Because ``arclength'' isn't formally defined, the strategy I'll take
is to make some reasonable assumptions that any notion of arclength should satisfy and work with them.

For this problem, and other problems which use geometric reasoning, it's important to note that the
Kock-Lawvere axiom can be stated in the following form: 

\begin{proposition}[Microstraightness] If $f\colon R\to R^n$ is any curve, $x\in R$, and
$d\in D$, then the portion of the curve from $f(x)$ to $f(x+d)$ is straight. 
\end{proposition}

Let $f\in R^R$ be any function, and let $s(x)$ be the arclength of the graph of $y = f(x)$ from 0 to $x$.
(That is, $s$ is the function which we would like to determine.) 

\begin{figure}
\begin{center}
\includegraphics{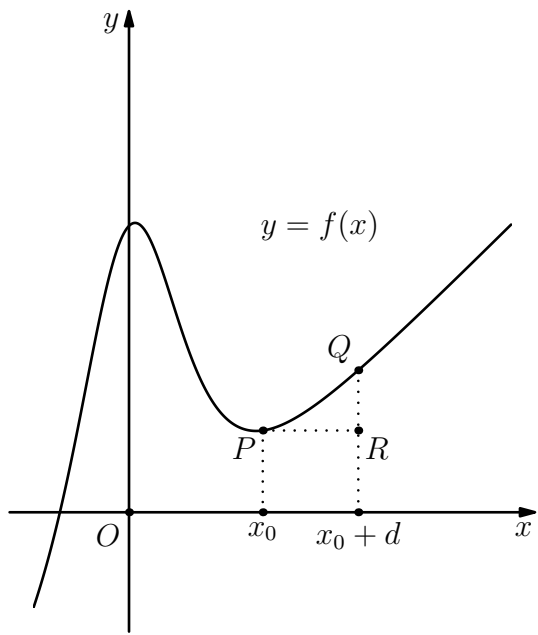}
\end{center}
\caption{Determining the Arclength of $y = f(x)$}
\label{arclength}
%
%
%
%
%
%
\end{figure}

Let $x_0 \in R$ and $d\in D$ be arbitrary and consider $s(x_0 + d) - s(x_0)$. It should be the length
of the segment of $y = f(x)$ from $x_0$ to $x_0 + d$, as in Figure~\ref{arclength}.

Because of microstraightness, we know that the part of the graph of $y=f(x)$ from $P$ to $Q$ is a straight line.  
Furthermore, it is the hypotenuse of a right triangle with legs $PR$ and $RQ$.  The length of $PR$ is $d$.

To determine the length of $RQ$: Note that the height of $P$ is $f(x)$, so the height of $R$ is $f(x)$. On
the other hand, the height of $Q$ is $f(x + d) = f(x) + f'(x)d$, so the length of $RQ$ is $f'(x)d$.

The hypotenuse of a right triangle with legs of length 1 and $f'(x)$ is $\sqrt{1 + f'(x)^2}$.  By scaling down,
we see that the length of $PQ$ is $d\sqrt{1+f'(x)^2}$.

So, we know that $s(x + d) - s(x)$ should be $d\sqrt{1 + f'(x)^2}$.  On the other hand, $s(x + d) - s(x) = 
ds'(x)$.  By microcancellation, we have that $s'(x) = \sqrt{1 + f'(x)^2}$.  Since $s(0) = 0$, we have
\begin{displaymath}
s(x) = \int_0^x \sqrt{1 + f'(t)^2}\,dt
\end{displaymath}

\begin{figure}
\begin{center}
\includegraphics{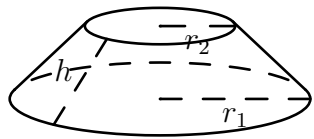}
\end{center}
\caption{A frustum of a cone}
\label{frustum}
%
%
%
%
%
%
\end{figure}

Several other formulas can be derived using precisely the same technique.  For example, suppose we want
to know the surface area of revolution of $y = f(x)$.  Furthermore, suppose we know that the surface area of
a frustum of a cone with radii $r_1$ and $r_2$ and slant height $h$ as in Figure~\ref{frustum} is $\pi(r_1 + r_2)h$. (See below to eliminate
this assumption.)

Then, let $A(x_0)$ be the surface area of revolution of $y = f(x)$ from $x = 0$ to $x = x_0$ about the $x$-axis.  As before,
consider $A(x_0 + d) - A(x_0)$ where $d$ is arbitrary.  This should be the surface area of the frustum
obtained by rotating $PQ$ about the $x$-axis.  The slant height is the length of $PQ$, which we determined
earlier was $(\sqrt{1 + f'(x)^2})d$.  The two radii are $f(x)$ and $f(x + d) = f(x) + f'(x)d$.  Therefore,
\begin{displaymath}
A(x_0 + d) - A(x_0) = \pi(f(x) + f(x) + f'(x)d)(\sqrt{1 + f'(x)^2})d
\end{displaymath}
which, multiplying out, becomes $d2\pi f(x)\sqrt{1 + f'(x)^2}$.  As before, $A(x_0 + d) - A(x_0)$ is
also equal to $A'(x_0)d$, so 
\begin{displaymath}
A(x) = 2\pi\int_0^x f(t) \sqrt{1 + f'(t)^2}\,dt
\end{displaymath}

\begin{exercise}Derive the 
formula for the volume of the solid of revolution of $y = f(x)$ about the $x$-axis.
\end{exercise}
\begin{exercise} Derive
the formula for the arclength of a curve $r = f(\theta)$ given in polar form.
\end{exercise}
\begin{exercise} Show that the (signed) area under the curve $y = f(x)$ from $x = a$ to
$x = b$ is $\int_a^b f(x)\,dx$.
\end{exercise}
\begin{exercise} Above we assumed that we knew the surface area of a frustum of a cone. Eliminate
this assumption by deriving the formula for the surface area of a cone (from which the 
formula for the surface area  of a frustum follows by an argument with similar triangles)
as follows: 

Fix a cone $C$ of slant height $h$ and radius
$r$.  The cone $C$ can be considered to be the graph of a function $y = mx$ from $x = 0$ to $x = r/m$
revolved a full $2\pi$ radians around the $x$-axis.  

Let $A(\theta)$ be the area of the surface formed by revolving the graph of $y = mx$ from $x = 0$
to $x = r/m$ only $\theta$ radians around the $x$-axis.

Using a method similar to that above, determine that $A(x) = (1/2)xrh$.  
This gives the surface area as
$A(2\pi) = \pi rh$.
\end{exercise}

\subsection{The Equation of a Catenary}

In the above section, essentially the same method was used again and again to solve different problems.  As
an example of a different way to apply SIA in single-variable calculus, in this section I'll outline
how the equation of a catenary may be derived in it. The full derivation is in~\cite{primer}.)

To do this, we'll need the existence of functions $\sin$, $\cos$, $\exp$ in $R^R$ satisfying
$\sin(0) = 0$, $\cos(0) = \exp(0) =1$, $\sin' = \cos$, $\cos' = -\sin$ and $\exp' = \exp$.  We get this
from the following set of axioms.

\newcommand{\R}{\mathbb{R}}

\begin{axiom}
For every $C^\infty$ function $f\colon \R^n \to \R^m$ (in the real world), we assume we have a function
$f\colon R^n \to R^m$ (in SIA).  Furthermore, for any true identity constructed out of such functions, composition, and partial differentiation operators, we may take the corresponding statement in SIA to be an axiom.  (``True'' means
true for the corresponding functions between cartesian products of $\R$ in the real world.)
\end{axiom}

(We can actually go further.  For every $C^\infty$ manifold $\mathbb{M}$ in the real world, we may assume
that there is a set $M$ in SIA, and for every $C^\infty$ function $f\colon \mathbb{M}\to\mathbb{N}$
we may assume that there is a function $f\colon M\to N$ in SIA, and we may assume that these functions
satisfy all identities true of them in the real world.  But I will not use these extra axioms in this article.)

\begin{figure}
\begin{center}
\includegraphics{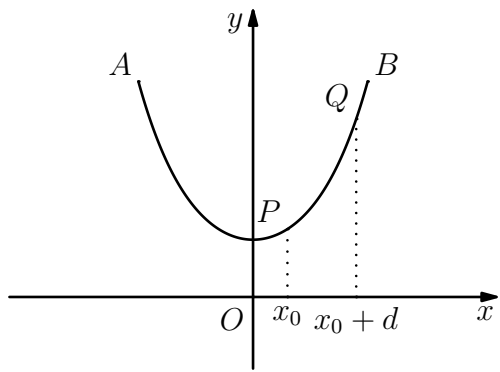}
\end{center} 
\caption{Determining the equation of a catenary}
\label{catenary}
\end{figure}

Suppose that we have a flexible rope of constant weight $w$ per unit length suspended from two points $A$ and
$B$ (see figure~\ref{catenary}).  We would like to find the function $f$ such that the graph of $y = f(x)$ is the curve that the rope makes. (We
will actually disregard the points $A$ and $B$ and consider $f$ to be defined on all of $R$.)

Let $T(x)$ be the tension in the rope at the point $(x,f(x))$.  (Recall that the tension at a point in
a rope in equilibrium is defined as follows: That point in the rope is being pulled by both sides of the
rope with some force.  Since the rope is in equilibrium, the magnitude of the two forces must be equal.
The tension is that common magnitude.)

Let $\phi(x)$ be the angle that the tangent to $f(x)$ makes with the positive $x$-axis. (That is,
$\phi(x)$ is defined so that $\sin\phi(x) = f'(x)\cos\phi(x)$). We suppose that we have chosen
the origin so that $\phi(0) = 0$.

Let $s(x)$ be the arclength of $f(x)$ from 0 to $x$.  

Let $x_0\in R$ and $d\in D$ be arbitrary.  Consider the segment of the rope from $P = (x_0,f(x_0))$
to $Q = (x_0 + d,f(x_0 + d))$.  This segment is in equilibrium under three forces:

1. A force of magnitude $T(x_0)$ with direction $\phi(x_0) + \pi$.

2. A force of magnitude $T(x_0 + d)$ with direction $\phi(x_0 + d)$.

3. A force of magnitude $w(s(x_0 + d) - s(x_0)) = ws'(x_0)d$ with direction $-\pi/2$.

\begin{exercise}
By resolving these forces horizontally and using microcancellation, show that the horizontal
component of the tension (that is, $T(x)\cos\phi(x)$) is constant.  Call the constant
tension $T_0$.
\end{exercise}
\begin{exercise}
By resolving these forces vertically and using microcancellation and the fact that $\phi(0) = 0$, show
that the vertical component of the tension (that is $T(x)\sin\phi(x)$) is $ws(x)$.
\end{exercise}
\begin{exercise}
By combining the previous two exercises and using the fact that $\sin\phi(x) = \cos\phi(x) f'(x)$
and $s'(x) = \sqrt{1 + f'(x)^2}$, show that $f$ satisfies the differential equation
$1 + (u')^2 = a^2(u'')^2$, where $a = T_0/w$.
\end{exercise}

Solving differential equations symbolically is the same in SIA as it is classically, since
no infinitesimals or limits are involved.  In this case, the answer turns out to be
\[f(x) = a\cosh\left(\frac{x}{a}\right) = \frac{a(e^{x/a} + e^{-x/a})}{2},\]
if we add the initial condition $f(0) = a$ to our previously assumed initial condition
$f'(0) = 0$.
\section{Multivariable Calculus}

\begin{definition}[Partial Derivatives]
Let $f(x,y)$ be a function from $R^2$ to $R$.  We define the partial derivative $\partial f/\partial x$
(also written $f_x$) as follows: Given $y$, let $g_y(x) = f(x,y)$.  Then $f_x(x_0,y_0)$ is 
defined to be $g'_{y_0}(x_0)$.  A similar definition is made for $f_y$, and for functions of more than
two variables.
\end{definition}

\begin{definition}[$D(n)$]
For $n\in \N$, let $D(n) = \{(d_1,\ldots, d_n)\in D^n \mid \forall i,j\, d_id_j = 0\}$.  Note that $D(1) = D$.
\end{definition}

The sets $D(n)$ play the role in multivariable calculus that $D$ played in singlevariable calculus.  For example,
we have the following.
\begin{proposition}
Let $f(x,y)$ be a function from $R^2$ to $R$.  Then, for all $(d_1, d_2) \in D(2)$, 
\begin{displaymath}
f(x_0 + d_1, y_0 + d_2) = f(x_0,y_0) + d_1 f_x(x_0,y_0) + d_2 f_y(x_0,y_0)
\end{displaymath}
and furthermore, $f_x(x_0,y_0)$ and $f_y(x_0,y_0)$ are unique with those properties.

The analogous statement is also true for functions of more than two variables.
\end{proposition}

We also have
\begin{proposition}[Extended Microcancellation]
Let $a_1, \ldots, a_n\in R$.  Suppose that for all $(d_1,\ldots,d_n) \in D(n)$, $\sum a_id_i = 0$.  Then each
$a_i$ equals 0.
\end{proposition}

\subsection{Stationary Points and Lagrange Multipliers}

There is an interesting substitute for the method of Lagrange multipliers in Smooth Infinitesimal Analysis.  To introduce it, I'll first discuss the concept of stationary points.

Suppose that we've forgotten what a stationary point and what a critical point is, and we need to redefine the
concept in Smooth Infinitesimal Analysis.  How should we do it?  We want a stationary point to be such that
every local maximum and local minimum is one.  A point $x$ gives rise to a local maximum $(x,f(x))$ of a 
single-variable function $f$ just in case there is some neighborhood of $x$ such that $f(x) \geq f(x_0)$
for all $x_0$ in that neighborhood.

However, in Smooth Infinitesimal Analysis, there is always a neighborhood of $x$ on which $f$ is \emph{linear}.
That means that for $x$ to be a local maximum, it must be \emph{constant} on some neighborhood.  Obviously,
the same is true if $x$ is a local minimum.  This suggests that we say that $f$ has a stationary point at
$x$ just in case $f(x) = f(x + d)$ for all $d\in D$.  

\begin{definition}[Stationary Point of a Single-Variable Function] Let $f\in R^R$ and $x\in R$.
We say that $f$ has a stationary point at $x$ if for all $d\in D$, $f(x + d) = f(x)$.
\end{definition}

Similarly, given a function $f(x,y)$ of two variables, and a point $(x_0,y_0)$, $f$ is linear on the set
$(x_0,y_0) + D(2)$.  This suggests the following definition. that we define $(x_0,y_0)$ to be a stationary point of $f$ just in
case $f(x_0,y_0) = f(x_0 + d_1, y_0 +d_2)$ for all $(d_1,d_2)\in D(2)$. 

\begin{definition}[Stationary Point of a Multivariable Function] Let $f\colon R^n\to R$.  We say that
$\bar{x}\in R^n$ is a stationary point of $f$ if for all $\bar{d}\in D(n)$, $f(\bar{x} + \bar{d}) = f(\bar{x})$.
\end{definition}

Now, suppose we want to maximize or minimize a function $f(x,y)$ subject to the constraint that it be
on some level surface $g(x,y) = k$, where $k$ is a constant.  Now, we should require of $(x_0,y_0)$
not that $f(x_0 + d_1, y_0 + d_2) = f(x_0,y_0)$ for all $(d_1, d_2) \in D(2)$, but only for those
$(d_1,d_2)\in D$ which keep $(x_0,y_0)$ on the same level surface of $g$; that is, those $(d_1,d_2)\in D(2)$
for which $g(x_0 + d_1, y_0 + d_2) = g(x_0,y_0)$.  I'll record this in a definition.

\begin{definition}[Constrained Stationary Point] Let $f$, $g\colon R^n\to R$. A point $\bar{x}\in R^n$ is
a stationary point of $f$ constrained by $g$ if for all $\bar{d}\in D^n$, if $g(\bar{x} + \bar{d}) = g(\bar{x})$
then $f(\bar{x} + \bar{d}) = f(\bar{x})$.
\end{definition}

I'll show how this definition leads immediately to a method of solving constrained extrema problems by doing an
example.

Suppose we want to find the radius and height of the cylindrical can (with top and bottom) of least surface
area that holds a volume of $k$ cubic centimeters.  The surface area is $f(r,h) = 2\pi r h + \pi r^2 + \pi r^2$,
and we are constrained by the volume, which is $g(r,h) = \pi r^2 h$.

We want to find those $(r,h)$ such that $f(r + d_1, h + d_2) = f(x_0,y_0)$ for all those
$(d_1,d_2)\in D(2)$ such that $g(r + d_1, h + d_2) = g(r,h)$. So, the first question is to
figure out which $(d_1,d_2)\in D(2)$ satisfy that property.

We have
\begin{displaymath}
g(r + d_1, h + d_2) = \pi (r + d_1)^2 (h + d_2)
\end{displaymath}
which is 
\begin{displaymath}
\pi(r^2 + 2rd_1)(h + d_2) = \pi(r^2h + 2rd_1h + r^2 d_2)
\end{displaymath}
If this is to equal $\pi r^2h$, then we must have $\pi(2rd_1h + r^2 d_2) = 0$,
so that $d_1 = -(r/(2h)) d_2$.

Now, we want to find an $(r,h)$ so that $f(r + d_1, h + d_2) = f(r,h)$ where $d_1 = -(r/(2h)) d_2$.

We have 
\begin{displaymath}
f(r + d_1, h + d_2) = 2\pi((r + d_1)(h + d_2) + (r + d_1)^2)
\end{displaymath}
which is
\begin{displaymath}
2\pi(rh + d_1h + d_2 r + r^2 + 2rd_1) = 2\pi(rh + r^2 + d_1(h + 2r) + d_2 r
\end{displaymath}
If this is to equal $2\pi(rh + r^2)$ then we must have $d_1(h + 2r) + d_2 r = 0$. Substiting
$d_1 = -(r/(2h)) d_2$, we get $(-(r/(2h))(h + 2r) + r)d_2 = 0$.  By microcancellation,
we have $-(r/(2h))(h + 2r) + r = 0$, from which it follows that $2r =h$.

\subsection{Stokes's Theorem}

It is interesting that not only can the theorems of vector calculus such as Green's theorem, Stokes's theorem,
and the Divergence theorem be stated and proved in Smooth Infinitesimal Analysis, but, just as in the classical
case, they are all special cases of a generalized Stokes's theorem.

In this section I will state Stokes's theorem.

\begin{definition} Given $x$, $y\in R$, we say that $x\leq y$ if $\neg(y < x)$.  We define
$[x,y]$ to be the set $\{z\in R\mid x \leq z \leq y \}$.
\end{definition}

\begin{definition} 
Let $C\colon [0,1]\to R^3$ be a curve, and $F = \langle M,N,P\rangle \colon R^3 \to R^3$ be a vector
field.  The line integral
\begin{displaymath}
\int_C F\cdot dr
\end{displaymath}
is defined to be $\int_0^1 F(C(t))\cdot C'(t)\,dt$.
\end{definition}

\begin{definition} Let $S = S(u,v)\colon [0,1]^2\to R^3$ be a surface, and $f\colon R^3 \to R$ be a function.  The surface
integral
\begin{displaymath}
\iint_S f\,d\sigma
\end{displaymath}
is defined to be $\int_0^1 \int_0^1 f(S(u,v))\cdot |S_u(u,v) \times S_v(u,v)|\,du\,dv$.
\end{definition}

This definition may be intuitively justified in the same manner that the arclength of a function was
derived in an earlier section.

\begin{definition} Let $S = S(u,v)\colon [0,1]^2\to R^3$ be a surface, and $F\colon R^3\to R^3$ be a vector
field.  The surface integral
\begin{displaymath}
\iint_S F\cdot n\,d\sigma
\end{displaymath}
is defined to be 
\[\iint_S F\cdot \left(\frac{S_u \times S_v}{|S_u \times S_v|}\right)\,d\sigma.\]

Note that this equals $\int_0^1 \int_0^1 F(S(u,v))\cdot (S_u(u,v) \times S_v(u,v))\,du\,dv$.
\end{definition}

We extend both definitions to cover formal $R$-linear combinations of curves and surfaces, and we define the boundary
$\partial S$ of a region $S$ to be the formal $R$-linear combination of curves $S(0,\cdot) + S(\cdot, 1) - S(1,\cdot) - S(\cdot,0)$.

The curl of a vector field $F = \langle M,N,P\rangle$ is defined as usual, and we can prove the usual Stokes's Theorem:
\begin{theorem} Let $S$ be a surface and $F$ a vector field.  Then
\begin{displaymath}
\iint_S \mathop{\mathrm{curl}} F \cdot n\,d\sigma = \int_{\partial S} F\cdot dr
\end{displaymath}
\end{theorem}

This theorem may be used to compute answers to standard multivariable calculus problems requiring Stokes's theorem
in the usual way.

\begin{exercise} State the divergence theorem in SIA.
\end{exercise}

\subsection{Generalized Stokes's Theorem}

The definitions in this section are directly from~\cite{moerdijk}.

\begin{definition}[Infinitesimal $n$-cubes] For $n\in \N$, and $S$ any set, an infinitesimal $n$-cube
in $S$ is some $(\bar{d},f)$ where $\bar{d}\in D^n$ and $f\colon D^n\to S$.
\end{definition}

Intuitively, an infinitesimal $n$-cube on a set $S$ is specified by saying how you want to map $D^n$ into
your set, and how far you want to go along each coordinate.

Note that an infinitesimal 0-cube is simply a point.

\begin{definition}[Infinitesimal $n$-chains] An infinitesimal $n$-chain is a formal $R$-linear combination
of infinitesimal $n$-cubes.
\end{definition}

\begin{definition}[Boundary of $n$-chains] Let $C$ be a 1-cube $(d,f)$.  The boundary $\partial C$ is
defined to be the 0-chain $f(d) - f(0)$, where this is a formal $R$-linear combination of 0-cubes identified as points.

Let $C$ be a 2-cube $((d_1,d_2),f)$.  The boundary $\partial C$ is defined to be the 1-chain
$(d_1,f(\cdot,0)) + (d_2,f(d_1,\cdot)) - (d_1,f(\cdot,d_2)) - (d_2,f(0,\cdot))$.

In general, if $C$ is an $n$-cube $(\bar{d},f)$, the boundary $\partial C$ is defined to be the
$n-1$-chain $\sum_{i = 1}^n \sum_{\alpha = 0,1} (-1)^{i + \alpha} ((d_1,\ldots,\hat{d_i},\ldots,d_n),(x_1,\ldots,x_n)
\mapsto f(x_1,\ldots, \alpha \cdot d_i,\ldots x_n))$.

The boundary map is extended to chains in the usual way.
\end{definition}

\begin{definition}[Differential Forms]
An $n$-form on a set $S$ is a mapping $\omega$ from the infinitesimal $n$-cubes on $S$ to $R$ satisfying

1. Homogeneity. Let $a\in R$, $1\leq i\leq n$, and $f\colon D^n\to S$.  Define $g\colon D^n\to S$ by
$g(\bar{d}) = f(d_1,\ldots,ad_i,\ldots,d_n)$.  Then for all $\bar{d}\in D^n$, $\omega((\bar{d},g)) = 
a\omega((\bar{d},f))$.

2. Alternation. Let $\sigma$ be a permutation of $\{1,2,\ldots,n\}$.  Then $\omega(\bar{d},\sigma f)
 = \mathrm{sgn}(\sigma)\cdot \omega(\sigma\bar{d},f)$, where $\sigma f = (x_1,\ldots, x_n)\mapsto
 f(x_{\sigma(1)},\ldots, x_{\sigma(n)})$ and $\sigma\bar{d} = (d_{\sigma(1)},\ldots, d_{\sigma(n)})$.
 
3. Degeneracy.  If $d_i = 0$, $\omega(\bar{d},f) = 0$.

We often write $\omega$ as
\[ C\mapsto \int_C\omega. \]
We extend $\omega$ to act on all $n$-chains in the usual way.
\end{definition}

These axioms intuitively say that $\omega$ is a reasonable way of assigning an oriented size to 
the infinitesimal $n$-cubes. 

The homogeneity condition says that if you double the length of one
side of an infinitesimal $n$-cube, you double its size.  

The alternation condition says that
if you swap the order of two coordinates in an infinitesimal $n$-cube, then you negate its oriented
size.  

The degeneracy condition says that if any side of the infinitesimal $n$-cube is of length 0, its oriented
size is of length 0.

By the Kock-Lawvere axiom, for all differential $n$-forms $\omega$, there is a unique map
$\tilde{\omega}\colon S^{D^n}\to R$ such that for all $\bar{d}\in D^n$
and $f\colon D^n\to S$ we have $\omega(\bar{d},f) = d_1\cdots d_n\cdot \tilde{\omega}(f)$.

\begin{definition}[Exterior Derivative]
The exterior derivative $d\omega$ of a differential $n$-form $\omega$ is an $n+1$-form defined by
\[\int_C d\omega = \int_{\partial C} \omega\]
for all infinitesimal $n+1$-cubes.
\end{definition} 

\begin{definition}[Finite $n$-cubes]
A finite $n$-cube in $S$ is a map $M$ from $[0,1]^n$ to $S$.

The boundary of a finite $n$-cube is defined in the same way that the boundary of an infinitesimal $n$-cube was
defined.
\end{definition}
In the above section, a curve was a finite 1-cube in $R^3$ and a surface was a finite 2-cube in $R^3$.
\begin{definition}[Integration of forms over finite cubes]
Let $\omega$ be an $n$-form on $S$ and $M$ a finite $n$-cube on $S$. Then
$\int_M \omega$ is defined to be
\[\int_0^1\cdots \int_0^1 \tilde{\omega}(\bar{d} \mapsto M(\bar{t} + \bar{d}))\,dt_1\ldots dt_n.\]
\end{definition}

Generalized Stokes's theorem (for finite $n$-cubes) is provable in SIA (see~\cite{moerdijk} for the proof).
\begin{theorem}[Generalized Stokes's Theorem]
Let $S$ be a set, $\omega$ an $n$-form on $S$, and $M$ a finite $n+1$-cube on $S$.  Then
\[ \int_{\partial M}\omega = \int_M d\omega \]
\end{theorem}

Let's see how this gives the Fundamental Theorem of Calculus.

Let $F\in R^R$ and let $f = F'$.  We would like to see how
\[ \int_0^1 f(t)\,dt = F(1) - F(0) \]
is a special case of Generalized Stokes's Theorem. (On the other hand, that it's \emph{true} is 
immediate from the way we defined integration.)

Let $\omega$ be the 0-form on $R$ defined by $\omega(x) = F(x)$.  (Recall that 0-cubes are identified with points.)

Then $d\omega$ is the 1-form which takes infinitesimal 1-cubes $(d,g)$ to $\int_{\partial(d,g)} \omega$.
We must show that for the finite 1-cube $[0,1]$, $\int_{[0,1]} d\omega = \int_0^1 f(t)\,dt$.

The boundary of $(d,g)$ is $g(d) - g(0)$ (as a formal linear combination, not as a subtraction in $R$).
Therefore, $\int_{\partial(d,g)} \omega = \omega(g(d)) - \omega(g(0)) = F(g(d)) - F(g(0))$.  
Since $g$ is a function from $D$ to $R$, there is a unique $a$ such that $g(d) = g(0) + ad$ for all $d\in D$.
Then $F(g(d)) - F(g(0)) = F(g(0) + ad) - F(g(0)) = F(g(0)) + (F'(g(0)))ad - F(g(0)) = f(g(0))ad$.  Therefore,
$\tilde{d\omega}(g) = f(g(0))a$, where $g(d) = g(0) + ad$ for all $d\in D$.

Therefore, $\int_{C}d\omega = \int_0^1 \tilde{d\omega}(d\mapsto d + t)\,dt = \int_0^1 f(t)\,dt$.

One can show in a similar manner that Stokes's theorem and the Divergence theorem are special cases
of Generalized Stokes's theorem, although the computations are significantly more arduous. 
%
%
%
%
%
%
%

\bibliography{sia}
\end{document}